\documentclass[12pt]{amsart}

\usepackage{eucal,url,amssymb, 
  enumerate,amscd}
\usepackage{accents}

\makeatletter
\def\widebar{\accentset{{\cc@style\underline{\mskip10mu}}}}
\def\Widebar{\accentset{{\cc@style\underline{\mskip8mu}}}}
\makeatother

\usepackage[pagebackref]{hyperref}

\theoremstyle{plain}
\newtheorem{thm}{Theorem}[section]
\newtheorem{theorem}{Theorem}[section]
\newtheorem{lemma}[thm]{Lemma}

\newtheorem{proposition}[thm]{Proposition}

\newtheorem{problem}{Problem}

\theoremstyle{definition}
\newtheorem{definition}[thm]{Definition}

\theoremstyle{remark}
\newtheorem{remark}[thm]{Remark}

\newcommand{\Lie}[1]{\operatorname{\textsl{#1}}}
\newcommand{\lie}[1]{{\operatorname{\mathfrak{#1}}}}
\DeclareMathOperator{\ad}{ad}
\newcommand{\GL}{\Lie{GL}}
\newcommand{\gl}{\lie{gl}}

\newcommand{\Un}{\Lie{U}}

\newcommand{\End}{\Lie{End}}

\newcommand{\Ad}{\Lie{Ad}}
\newcommand{\Ann}{\textup{Ann}}

\newcommand{\CC}{\mathbb{C}}

\newcommand{\oper}[2]{\newcommand{#1}{\mathop{\mathrm{#2}}\nolimits} }
\oper{\DGA}{DGA}

\def\qed{\rule{2.3mm}{2.3mm}}
\newcommand{\bproof}{\noindent{\it Proof: }}
\newcommand{\eproof}{\hfill \qed \vspace{0.2in}}

\newcommand{\inp}[2]{\left< #1, #2\right>}

\DeclareMathOperator{\dbar}{\bar\partial}
\DeclareMathOperator{\bw}{\wedge}

\title[Weak Mirror Symmetry]%
{Weak Mirror Symmetry of Lie Algebras}

\date{April 30, 2008}

\author{R. Cleyton}
\address[Cleyton]{Institute f\"ur Mathematik\\
Humboldt Universit\"at zu Berlin\\Unter den Linden 6\\
10099 Berlin\\
Germany}
\email{cleyton@mathematik.hu-berlin.de}

\author{J. Lauret}
\address[Lauret]{FaMAF and CIEM\\ Universidad Nacional de C\'ordoba\\ 5000 C\'ordoba, Argentina}
\email{lauret@mate.uncor.edu}

\author{Y. S. Poon}
\address[Poon]{Department of Mathematics\\
University of California, Riverside\\ CA 92521\\
USA} \email{ypoon@ucr.edu}

\begin{document}

\maketitle

\section{Introduction}

It is well known that deformation theory of geometric objects such as
complex structures and symplectic structures are dictated by a
differential Gerstenhaber algebra (a.k.a. DGA) and the associated
cohomology theory \cite{GM}, \cite{Zhou}.  Therefore, DGA plays a key
role in mirror symmetry \cite{BarannikovKontsevich} \cite{Fukaya}. In
developing the algebraic aspects of mirror symmetry, Merkulov proposes
the concept of \it weak mirror symmetry \rm \cite{Mer-note}.  If $M$
is a manifold with a complex structure $J$ and $M^\vee$ is another
manifold with a symplectic structure $\omega$, then $(M, J)$ and
$(M^\vee, \omega)$ form a weak mirror pair if the associated
differential Gerstenhaber algebras $\DGA(M, J)$ and $\DGA(M^\vee,
\omega)$ are quasi-isomorphic.  The overall goal of this project is to
construct all mirror pairs when the manifolds $M$ and $M^\vee$ are
solvmanifolds, i.e. homogeneous spaces of simply-connected connected
solvable Lie groups, $J$ is an invariant complex structure and
$\omega$ is an invariant symplectic structure.

In the SYZ-conjecture, one considers the geometry of special
Lagrangian fibrations in a Calabi-Yau manifold $L\hookrightarrow
M\rightarrow B$ with $L$ being a real three-dimensional torus. The
mirror image is presumably a (new) Calabi-Yau manifold $M^\vee$ with
the special Lagrangian fibrations $L^*\hookrightarrow
M^\vee\rightarrow B$ where $L^*$ is the dual torus of $L$ \cite{SYZ}.
One way to adapt the structure of a Lie group $H$ to resemble this
situation is by insisting that the Lie algebra $\lie h$ is a
semi-direct sum of a subalgebra $\lie g$ by an abelian ideal $V$.  The
group $H$ is then a semi-direct product, namely the product of the
group $G$ corresponding to $\lie g$ with $V$.  By restricting our
attention to invariant structures on homogenous spaces of such groups
the geometry of the fibration of $H$ over $G$ becomes encoded in the
corresponding objects on the Lie algebra $\lie h$ of $H$.  We shall
speak somewhat sloppily of $\lie g$ as the base and of $V$ as the
fiber of the fibration.

Forgetting about the SYZ-conjecture, semi-direct products are still
natural objects to study in connection with ``weak mirror symmetry''.
This is so since the direct sum of the bundle of type $(1,0)$ vectors
$T^{(1,0)}$ and bundle of $(0,1)$-forms $T^{*(0,1)}$ on a complex
manifold carries a natural Lie bracket (Schouten) such that
$T^{(1,0)}$ is a sub-algebra and $T^{*(0,1)}$ is an abelian ideal. It
is the associated exterior algebra of this semi-direct sum and
associated $\dbar$ complex that eventually controls the deformations
of the complex structure.

It is well-known that a symplectic structure defines a flat
torsion-free connection on a Lie algebra $\lie g$~\cite{Chu}.  As we
shall see, a flat torsion-free connection on a Lie algebra $\lie g$
also defines a symplectic form $\omega$, not on $\lie g$ but on a
semi-direct product $\lie h^\vee$ of $\lie g$ with its dual $\lie
g^*$, see also~\cite{MR994379}.  This symplectic form is defined such
that it is Lagrangian with respect to both the base $\lie g$ and the
fiber $V^*=\lie g^*$ and so we call the pair $(\lie h^\vee,\omega)$ a
\emph{Lagrangian semi-direct product}.  A torsion-free flat connection
on $\lie g$ also defines a \emph{totally real semi-direct product}
$(\lie h,J)$ where $J$ is a complex structure on a semi-direct sum of
$\lie g$ with itself.  Such complex structures are particular cases of
\emph{complex product structures}, see~\cite{MR2130906}.  The
observation that complex structures and symplectic structures on
certain semi-direct products both are related to the notion
torsion-free flat connections may be found in~\cite{MR2104679}.
Torsion-free flat connections are also known as affine structures and
as such already has widespread application in the study of mirror
symmetry, see for instance~\cite{AP, BB-1, gross-2007}.
Left-invariant torsion-free flat connections on Lie groups are
equivalent to (Lie compatible) left-symmetric algebras.  Much is known
about left-symmetric algebras.  In particular existence problems have
been examined and it is known that no left-symmetric structure exists
on semi-simple algebras. Also, certain nilpotent algebras of dimension
greater than nine have no left-symmetric structures, see
e.g.~\cite{MR1608273,MR1316552,MR1411303,MR1612015}.

Therefore we confine the scope of our present paper to deal with
solvable spaces, i.e. the base $\lie g$ and hence the total space
$\lie h$ are solvable algebras.  Recent advance in resolving the
Benson-Gordon conjecture means that when we insist that $\lie h$
should carry a K\"ahler structure then $\lie h$ is flat and therefore
of a very special solvable type \cite{BG-2,B-Cortes,Hase-2,MR0425012}.
For a nilpotent algebra more is true - it is K\"ahler only if it is
abelian \cite{BG-1, Hase-1}.  Therefore an invariant symplectic form
$\omega$ on a non-abelian nilpotent algebra is of type $(1,1)$ with
respect a complex structure $J$ if and only if $\omega$ and $J$
determine a non-definite metric.  We call such a pair a
\emph{pseudo-K\"ahlerian geometry}.

Invariant complex structures, their Dolbeault cohomology and
pseudo-K\"ahlerian geometry on nilmanifolds have been a subject of
much investigation in recent years, especially when the complex
dimension is equal to three, see
\cite{CF,CFP,CFGU,CFGAU,CFU,Rolle,Salamon}. In particular, with the
recent advance in understanding the cohomology theory on nilmanifolds
\cite{Rolle}, our computation and results on nilpotent algebras in
this paper could be used to provide a full description of the
differential Gerstenhaber algebras of any invariant complex structures
on nilmanifolds in all dimension.

We have organized ourselves as follows.  In the next section, we
briefly review the construction of differential Gerstenhaber algebras
for complex and symplectic structures, the definition of semi-direct
products and establish notations for subsequent computation. In
Section \ref{sec:geometry}, we study the complex and symplectic
geometry on semi-direct products.  We first establish a correspondence
between a totally real semi-direct product and flat torsion-free
connections on the base in Proposition \ref{prop:complex and
  connection}.  The analogous result for Lagrangian semi-direct
products is obtained in Proposition \ref{prop:symplectic and
  connection}.

In Section \ref{sec:dual}, we develop the concept of dual semi-direct
product and use it as the candidate of ``mirror'' space.  In Section
\ref{sec:special lag}, we demonstrate a construction of a special
Lagrangian structure on the dual semi-direct product whenever a
special Lagrangian structure on a semi-direct product is given. This
is the result of Proposition \ref{prop:special}. After a brief revisit
to the subject on flat connections in Section \ref{sec:flat again}, we
prove the first main theorem (Theorem \ref{thm:iso}), which states
that the differential Gerstenhaber algebra on a totally real
semi-direct product is isomorphic to the differential Gerstenhaber
algebra of the Lagrangian dual semi-direct product as constructed in
Proposition \ref{prop:complex and connection}.

In Section \ref{sec:examples}, we exhibit with some examples of
K\"ahlerian solvable algebras and their mirror partners. In Section
\ref{sec:nil}, we focus on nilpotent algebras of dimension-four and
dimension-six. The work on the four-dimensional case is a brief review
of past results \cite{Poon}.  Our first step in addressing the issue
of finding mirror pairs of special Lagrangian nilpotent algebras in
dimension-six begins in Section \ref{sec:algebraic}. In this section,
we determine nilpotent algebras admitting a semi-direct product
structures, and then identify their dual semi-direct product space in
Table (\ref{tab:algebraic mirror}).  In Section \ref{sec:geometric},
using results in literature we identify the semi-direct product
structures which potentially could admit totally real complex
structures or Lagrangian symplectic structures. The result is in Table
(\ref{tab:geometric mirror}). We finish this paper by giving examples
of special Lagrangian pseudo-K\"ahler structures on every algebras in
Table (\ref{tab:geometric mirror}), and identifying their mirror
structures.  This is the content of Theorem \ref{thm:mirror}.

\section{Preliminaries}

We first recall two well known constructions of differential
Gerstenhaber algebras (DGA) \cite{CP, Mer-note, Zhou}.  After a
motivation due to weak mirror symmetry, we recall the definition of
semi-direct product of Lie algebras.

\subsection{DGA of a complex structure}
Suppose $J$ is an integrable complex structure on $\mathfrak{h}$.
i.e. $J$ is an endomorphism of $\mathfrak{h}$ such that $J\circ J=-1$
and
\begin{equation} \label{integrable} [ {x} \bullet {y} ]+J[ {Jx}\bullet
  {y} ]+J[ {x} \bullet {Jy} ]-[ {Jx} \bullet {Jy} ]=0.
\end{equation}
Then the $\pm i$ eigenspaces $\mathfrak{h}^{(1,0)}$ and
$\mathfrak{h}^{(0,1)} $ are complex Lie subalgebras of the
complexified algebra $\mathfrak{h}_\CC$. Let $\mathfrak{f}$ be the
exterior algebra generated by $\mathfrak{h}
^{(1,0)}\oplus\mathfrak{h}^{*(0,1)}$, i.e.
\begin{equation}
  \mathfrak{f}^n:=\wedge^n(\mathfrak{h}^{(1,0)}\oplus\mathfrak{h}^{*(0,1)}),
  \quad \mbox{ and } \quad \mathfrak{f}=\oplus_n\mathfrak{f}^n.
\end{equation}
The integrability condition in (\ref{integrable}) implies that
$\mathfrak{f}^1$ is closed under the \emph{Schouten bracket}
\begin{equation}
  [ x+\alpha \bullet y+\beta ]:=[x, y]+\iota_xd\beta-\iota_yd\alpha.
\end{equation}
Note that working on a Lie algebra, the Schouten bracket coincides
with \emph{Courant bracket} in Lie algebroid theory.

A similar construction holds for the conjugate
${\overline{\mathfrak{f}}}$, generated by $ \mathfrak{h}^{(0,1)}
\oplus \mathfrak{h}^{*(1,0)}.  $

Let $d$ be the Chevalley-Eilenberg (C-E) differential $d$ for the Lie
algebra $\lie h, [-\bullet -]$. Then $(\wedge \mathfrak{h}^*, d)$ is a
differential graded algebra. Similarly, let $\dbar$ be the C-E
differential for the complex Lie algebra
${\overline{\mathfrak{f}}}^1$.  Note that the natural pairing on
$(\lie h\oplus\lie h^*)\otimes \mathbb C$ induces a complex linear
isomorphism $ ({\overline{\mathfrak{f}}}^1)^*\cong
{{\mathfrak{f}}}^1$.  Therefore, the C-E differential of the Lie
algebra ${\overline{\mathfrak{f}}}^1$ is a map from $\mathfrak{f}^1$
to $\mathfrak{f}^2$.  Denote this operator by $\overline\partial$. It
turns out that $(\mathfrak{f}, [ {-} \bullet {-} ], \wedge,
\overline\partial)$ form a differential Gerstenhaber algebra which we
denote by $\DGA (\lie h, J)$.  The same construction shows that
$({\overline{\mathfrak{f}}}, [ {-} \bullet {-} ], \wedge, \partial)$
is a differential Gerstenhaber algebra, conjugate linearly isomorphic
to $\DGA(\lie h, J)$.  The above construction could be carried out
similarly on a manifold with a complex structure.

\subsection{DGA of a symplectic structure}
Let $\mathfrak{k}$ be a Lie algebra over $\mathbb{R}$. Suppose that
$\omega$ is a symplectic form on $\lie k$. Then the contraction with
$\omega$, $\omega:\mathfrak{k}\to\mathfrak{k}^*$ is a real
non-degenerate linear map. Define a bracket $[ {-} \bullet {-}
]_\omega$ on $\mathfrak{k}^*$ by
\begin{equation} [ {\alpha} \bullet {\beta} ]_\omega:=\omega[
  {\omega^{-1}\alpha} \bullet {\omega^{-1}\beta} ].
\end{equation}
It is a tautology that $(\mathfrak{k}^*, [ {-} \bullet {-} ]_\omega)$
becomes a Lie algebra, with the map $\omega$ understood as a Lie
algebra homomorphism.

In addition, the exterior algebra of the dual $\mathfrak{h}^*$ with
the C-E differential $d$ for the Lie algebra $\lie k$ is a
differential graded Lie algebra. In fact, $(\wedge^\bullet{\lie k}^*,
[-\bullet -]_\omega, \wedge, d)$ is a differential Gerstenhaber
algebra over $\mathbb{R}$.  After complexification we denote this by
$\DGA ({\lie k}, \omega)$.

\subsection{Quasi-isomorphisms and isomorphisms}
\begin{definition} {\rm \cite{Mer-note}} The Lie algebra $\lie h$ with
  an integrable complex structure $J$ and the Lie algebra $\lie k$
  with a symplectic structure $\omega$ form a weak mirror pair if the
  differential Gerstenhaber algebras $\DGA(\lie h, J)$ and $\DGA(\lie
  k, \omega)$ are quasi-isomorphic.
\end{definition}

Suppose that $\phi: \DGA(\lie k, \omega) \to \DGA(\lie h, J)$ is a
quasi-isomorphism. Since the concerned DGAs are exterior algebras
generated by finite dimensional Lie algebras, it is natural to examine
the property of the Lie algebra homomorphism on the degree-one
elements.
\begin{equation}\label{phi 1}
  \phi: \lie k^*_\CC \to \lie
  f^1=\mathfrak{h}^{(1,0)}\oplus\mathfrak{h}^{*(0,1)}.
\end{equation}
In particular, if the restriction of $\phi$ to $\lie k^*_\CC$ is an
isomorphism, it induces an isomorphism from $\DGA(\lie h, J)$ to
$\DGA(\lie k, \omega)$.  It turns out that for a special class of
algebras, this is the only situation when quasi-isomorphism occurs.

\begin{proposition}{\rm \cite{CP}}
  Suppose that $\lie h$ and $\lie k$ are finite dimensional nilpotent
  Lie algebras of the same dimension, $J$ is an integrable complex
  structure on $\lie h$ and $\omega$ is a symplectic form on $\lie k$.
  Then a homomorphism $\phi$ from $\DGA ({\lie h}, J)$ to $\DGA ({\lie
    k}, \omega)$ is a quasi-isomorphism if and only if it is an
  isomorphism.
\end{proposition}

This proposition provides a large class of Lie algebras to work on.
So in this paper we focus our attention on a restricted type of weak
mirror pairs.  Namely, we seek a pair such that the map $\phi$ in
(\ref{phi 1}) is an isomorphism on the degree-one level. Since the Lie
algebra $\lie k^*$ is tautologically isomorphic to $\lie k$ via
$\omega$, we concern ourselves with the non-degeneracy of the map
\begin{equation}
  \phi\circ \omega: \lie k_\CC\to \lie k^*_\CC \to \lie
  f^1=\mathfrak{h}^{(1,0)}\oplus\mathfrak{h}^{*(0,1)}.
\end{equation}
When this is an isomorphism one immediately obtains conditions on the
structure of $\lie k$. The reason for this is that with respect to the
Schouten bracket, $\mathfrak{h}^{(1,0)}$ is a subalgebra of $\lie f^1$
and $\mathfrak{h}^{*(0,1)}$ is an abelian ideal. In addition, $\dim
\mathfrak{h}^{(1,0)}=\dim \mathfrak{h}^{*(0,1)}$. In other words,
$\lie f^1$ and $\lie k$ are semi-direct products of a very particular
form.

\subsection{Semi-direct products}

Let $\lie g$ be a Lie algebra, $V$ a vector space and let
$\rho\colon\lie g\to\End(V)$ be a representation.  On the vector space
$\lie g\oplus V$, define
\begin{equation}
  \label{eq:1}
  [x+u,y+v]_\rho:=[x,y]+ \rho(x)v-\rho(y)u,
\end{equation}
where $x,y$ are in $\lie g$ and $u,v$ are in $V$.  Then this
determines a Lie bracket on $\lie g\oplus V$.  This structure is a
particular case of a semi-direct product.  As a Lie algebra it is
denoted by $\lie h=\lie h(\lie g,\rho)=\lie g\ltimes_\rho V$. By
construction, $V$ is an abelian ideal and $\lie g$ is a complementary
subalgebra.  One may also consider the semi-direct product as an
extension of the algebra $\lie g$ by a vector space $V$.

Note that if $H$ is a simply connected, connected Lie group of $\lie
h$ and $G$ the connected subgroup of $H$ with $\lie g\leqslant \lie h$
as above then $H/G\simeq \mathbb R^{\dim V}$ as a flat symmetric
space.

Conversely, suppose $V$ is an abelian ideal of a Lie algebra $\lie h$
and $\lie g$ a complementary subalgebra.  The adjoint action of $\lie
g$ on $V$ then gives $\rho\colon \lie g\to\End(V)$.

In this paper, we are solely interested in the situation when $\dim
\lie g=\dim V$. In particular, when the vector space $V$ is regarded
as the underlying space of the real algebra $\lie g$ or its dual $\lie
g^*$, interesting geometry and other phenomena arise through the
representations of $\lie g$ as described next.

\section{Geometry on semi-direct products}\label{sec:geometry}

A \emph{left-invariant connection} on a Lie group $G$ is an affine
connection $\nabla$ such that $\nabla_{X_g}Y_g=(L_g)_*(\nabla_{x}y)$,
where $X_g=(L_g)_*x,~Y_g=(L_g)_*y$, i.e. such that covariant
differentiation of left-invariant vector fields give rise to
left-invariant vector fields.  These are in one-to-one correspondence
with linear maps $\gamma\colon\lie g\to\End(\lie g)$ through
$\nabla_xy=\gamma(x)y$.  The torsion of $\gamma$ is
\begin{gather*}
  T^\gamma(x,y):=[x,y]-\gamma(x)y+\gamma(y)x
\end{gather*}
and its curvature $R^\gamma$
\begin{gather*}
  R^{\gamma}(x,y) := \gamma([x,y]) - \gamma(x)\gamma(y) +
  \gamma(y)\gamma(x).
\end{gather*}
Since all connections considered here are left-invariant, linear maps
$\gamma\colon \lie g\to\End(\lie g)$ are referred to as connections on
$\lie g$ and say that $\gamma$ is \emph{flat} if $R^\gamma=0$ and
\emph{torsion-free} if $T^\gamma=0$.

\subsection{Totally real semi-direct products}

A complex structure on a real vector space $W$ is an endomorphism $J$
such that $J^2=-1$.  If $V$ is a subspace of $W$ such that $V'=JV$
satisfies $V\oplus V'=W$ then we say that $J$ is \emph{totally real
  with respect to $V$}.  Given a totally real $J$ any $w$ in $W$ may
be written uniquely as $w=x+Jy$ for $x,y$ in $V$.  So $W\cong V\oplus
V$ and $J$ may be viewed as the standard complex structure $J_0$ on
$V\oplus V$ given by $(x,y)\mapsto (-y,x)$.  Since a basis for $W$
always may be chosen so that $Je_{2i-1}=e_{2i}$ any $J$ is totally
real with respect to some $V$.

If $\lie g=(W,[\cdot,\cdot])$ is a Lie algebra we say that $J$ is
\emph{integrable} if the Nijenhuis tensor
\begin{equation*}
  N_J(x,y):=[x,y]-[Jx,Jy]+J([x,Jy]+[Jx,y])
\end{equation*}
is zero for all $x$ and $y$ in $W$.  If $J$ is totally real with
respect to $V$ then $J$ is integrable if and only if $N_J(x,y)=0$ for
all $x,y\in V$. This follows by the identity $N_J(x,y) = JN_J(x,Jy)$
valid for all $x,y\in W$.

\begin{definition} Suppose that $\lie h(\lie g,\rho)=\lie
  g\ltimes_\rho V$ is a semi-direct product Lie algebra. A complex
  structure on $\lie g\ltimes_\rho V$ is totally real if $J$ is
  totally real with respect to $\lie g$ and $J\lie g=V$.
\end{definition}

Since $V$ is an abelian ideal, $[Jx,Jy] = 0$ for all $x,y$ in $\lie g$
and so the Nijenhuis tensor vanishes precisely when
\begin{equation}
  \label{eq:5}
  [x,y]+J\rho(x)Jy-J\rho(y)Jx=0.
\end{equation}
for all $x,y\in\lie g$.   This has the significance that
\begin{equation}
\gamma(x)y:=-J\rho(x)Jy
\end{equation}
defines a torsion-free connection on $\lie g$.  This is flat since
\begin{eqnarray*}
  &&  \gamma([x,y]) - \gamma(x)\gamma(y) +
  \gamma(y)\gamma(x) \\
  &=& -J\rho([x,y])J + J\rho(x)\rho(y)J -
  J\rho(y)\rho(x)J = 0.
\end{eqnarray*}

On the other hand, take a flat, torsion-free connection $\gamma$ on
$\lie g$.  Then the totally real complex structure $J$ on $\lie
h:=\lie g\ltimes_{\gamma}\lie g$ defined by $J(x,y)=(-y,x)$ becomes
integrable with respect to $[\cdot,\cdot]_\gamma$ by virtue of
\begin{align*}
  N_J((x,0),(y,0)) & = [(x,0),(y,0)]_\gamma +
  J[(x,0),(0,y)]_\gamma - J[(y,0),(0,x)]_\gamma\\
  & = [(x,0),(y,0)]_\gamma +
  J((0,\gamma(x)y) - (0,\gamma(y)x))\\
  &= ([x,y] - \gamma(x)y + \gamma(y)x, 0) = (0,0).
\end{align*}
This proves our first Proposition, which is at least implicitly
contained in \cite{MR2130906}.

\begin{proposition}\label{prop:complex and connection}
  There is a one-to-one correspondence between flat torsion-free
  connections on $\lie g$ and totally real integrable complex
  structures on semi-direct products $\lie g\ltimes_\rho V$.
\end{proposition}

\subsection{Lagrangian semi-direct products} Let $\omega$ be a
two-form on a vector space $W$, i.e. $\omega\in \Lambda^2 W^*$. We may
also view $\omega$ as a linear map $\omega\colon W\to W^*$ such that
$\omega^*=-\omega$ through the identification $(W^*)^*=W$. Let
$2m=\dim W$.  Then $\omega$ is \emph{non-degenerate} if
$\omega^{m}\not=0$. Equivalently, $\omega\colon W\to W^*$ is
invertible.

A subspace $V$ of $W$ is \emph{isotropic} if $\omega(V,V)=0$.
Equivalently $\omega(V)\subset \Ann(V)\subset W^*$.  If $\dim
V=\frac12 W$, then an isotropic $V$ is called a \emph{Lagrangian}
subspace. In this case $\omega(V)=\Ann (V)$.  A splitting $W=V\oplus
V'$ of a vector space $W$ into a direct sum is called Lagrangian with
respect to $\omega$ if both $V$ and $V'$ are Lagrangian with respect
to $\omega$.  Two vector spaces $W_1$ and $W_2$ with non-degenerate
two-form $\omega_1$ and $\omega_2$ are said to be isomorphic if a
linear isomorphism $f\colon W_1\to W_2$ exists such that
$f^*\omega_2=\omega_1$.

Let $V$ be a vector space.  Then $V\oplus V^*$ carries a two-form
$\omega$ given by the canonical pairing $\omega(\alpha,x)=\alpha(x)$
for $\alpha\in V^*$ and $x\in V$ and such that the given splitting is
Lagrangian.

\begin{lemma}
  Suppose $W$ is a vector space and $\omega$ is a non-degenerate
  two-form on $W$.  Then $W$ admits a Lagrangian splitting $W=V\oplus
  V'$ if and only if $(W,\omega)$ is isomorphic to $(V\oplus V^*,
  \langle\cdot,\cdot\rangle)$ (where $V$ may be taken to be $\mathbb
  R^{\frac12\dim W}$.
\end{lemma}
\bproof Clearly, if $f:W\to V\oplus V^*$ is an isomorphism then the
splitting given $U:=f^{-1}(V),~U':=f^{-1}(V^*)$ is Lagrangian with
respect to $\omega$.  If, on the other hand $U\oplus U'$ is a
Lagrangian splitting of $W$ then $\Ann(V')=\omega(V')$.  Furthermore
$\Ann(V')$ is canonically isomorph to $V^*$.  Therefore $W\cong
V\oplus V^*$ by the map $f(x+x'):=x+\omega(x')$.  Moreover,
\begin{align*}
  \inp{f(x+x')}{f(y+y')}&=\inp{x+\omega(x')}{y+\omega(y')}\\
  &=\omega(x',y)+\omega(x,y')\\
  &=\omega(x+x',y+y').
\end{align*}
\eproof

Suppose that $\lie g=(W,[\cdot,\cdot])$ is a Lie algebra.  Then the
derivative of a two-form $\omega$ with respect to the
Chevalley-Eilenberg differential is
\begin{align*}
  (d\omega)(x,y,z) &= -\left( \omega([x,y],z) + \omega([y,z],x) +
    \omega([z,x],y)\right)\\
  &= - \left( \omega([x,y])z - \omega(x)(\ad(y)z) +
    \omega(y)(\ad(x)z)\right).
\end{align*}
So $\omega$ is closed if and only if for all $x,y$:
\begin{equation}\label{eq:2}
  \omega([x,y]) = \ad^*(x)(\omega(y)) - \ad^*(y)(\omega(x))
\end{equation}
where $(\ad^*(x)\alpha)(y)=-\alpha([x,y])$.

\begin{definition} Suppose that a Lie algebra $\lie h$ is a
  semi-direct product $\lie h=\lie g\ltimes_\rho V$. It is said to be
  Lagrangian with respect to a non-degenerate 2-form if the subalgebra
  $\lie g$ and the abelian ideal $V$ are both Lagrangian with respect
  to $\omega$.
\end{definition}

When $\lie g\ltimes_\rho V$ is Lagrangian, there is a canonical
isomorphism $\omega(V)=\Ann(V)\cong \lie g^*$.  Similarly,
$\omega(\lie g)=\Ann(\lie g)\cong V^*$. Define
\begin{equation}
\rho^*\colon \lie g\to \End(V^*) \quad \mbox{ by } \quad
(\rho^*(x)\alpha)(u) = -\alpha(\rho(x)u).
\end{equation}
Then
\[
(\rho^*(x)\omega(y))(u) = -\omega(y,\rho(x)u)=-\omega(y,[x,u]_\rho) =
(\ad^*_\rho(x)\omega(y))(u).
\]
Comparing to equation~\eqref{eq:2} it is now clear that $\omega$ is
closed if and only if
\begin{equation}\label{closedness}
  \omega([x,y]) = \rho^*(x)(\omega(y)) - \rho^*(y)(\omega(x))
\end{equation}
for all $x,y\in\lie g$.  The story now repeats itself. Define
\begin{equation}
  \gamma(x)y:=\omega^{-1}(\rho^*(x)\omega(y)).
\end{equation}
This defines a flat torsion-free connection on $\lie g$ since
\begin{eqnarray*}
  &&  R^\gamma(x,y)z \\
  &=& \omega^{-1}(\rho^*([x,y])\omega(z)) \\
  && -
  \omega^{-1}(\rho^*(x)\omega(\omega^{-1}(\rho^*(y)\omega(z)))) +
  \omega^{-1}(\rho^*(y)\omega(\omega^{-1}(\rho^*(x)\omega(z))))\\
  &=& \omega^{-1}(\rho^*([x,y])\omega(z)) -
  \omega^{-1}(\rho^*(x)\rho^*(y)\omega(z)) +
  \omega^{-1}(\rho^*(y)\rho^*(x)\omega(z))\\
  &=&0.
\end{eqnarray*}
Conversely, take a flat torsion-free connection $\gamma$ on $\lie g$.
Let $\omega$ be the standard skew pairing on $\lie g\oplus\lie g^*$:
\[
\omega(x+u,y+v) = u(y) - v(x).
\]
Define the bracket on $\lie g\oplus\lie g^*$ as the semi-direct
product by representation $\gamma^*$.  Then the semi-direct product is
Lagrangian with respect to $\omega$. It follows that if $x,y,z$ are in
$\lie g$ and $u,v,w$ are in $\lie g^*$,
\begin{equation}\label{vanishing components}
  d\omega(x,y,z)=d\omega(u,v,w)=d\omega(x,u,v)=0.
\end{equation}
Moreover,
\begin{align*}
  (d\omega)(x,y,u) &= -\left(\omega([x,y],u) + \omega(\gamma^*(y) u,x)
    -
    \omega(\gamma^*(x)u,y)\right)\\
  &= \left(u([x,y]) - (\gamma^*(y)u)(x) +(\gamma^*(x)u)(y)\right)\\
  &=  u([x,y]) +  u(\gamma(y)x) - u(\gamma(x)y)\\
  &= u([x,y]-\gamma(x)y+\gamma(y)x) = 0.
\end{align*}
This gives us the following result (which may also be found
in~\cite{MR994379}).
\begin{proposition}\label{prop:symplectic and connection}
  There is a one-to-one correspondence between flat torsion-free
  connections $\rho$ on Lie algebras $\lie g$ and Lagrangian
  semi-direct products $\lie g\ltimes_\rho V$.
\end{proposition}

\subsection{From complex structure to two-form, and back}
\label{sec:dual}

Proposition \ref{prop:complex and connection} and Proposition
\ref{prop:symplectic and connection} of the preceding sections yield a
one-to-one correspondence between certain integrable complex
structures and certain symplectic forms going via flat, torsion-free
connections. In this section, we construct a direct relation between
two-forms and complex structures.

Suppose $W=V\oplus JV$, i.e. $J$ is a complex structure on $W$ so that
$V$ and $JV$ are totally real. On $W^\vee:=V\oplus (JV)^*$, define
\begin{equation}\label{omega-j}
  \omega_J(x+u,y+v) := v(Jx) - u(Jy),
\end{equation}
where $x,y$ are in $V$ and $u,v$ are in $(JV)^*$. Then $\omega_J$ is
non-degenerate on $W^\vee$ with both $V$ and $(JV)^*$ being
Lagrangian.

Conversely suppose $\omega$ is a non-degenerate 2-form on $W=V\oplus
V'$ with both $V$ and $ V'$ being Lagrangian. Write
$V'=\omega^{-1}(V^*)$ and set
\begin{equation}\label{j-omega}
  J_\omega(x+u)= - \omega^{-1}(u) + \omega(x)
\end{equation}
for all $x+u$ in $V\oplus V^*$.  Clearly, both $V$ and $V^*$ are
totally real with respect to $J_\omega$.

When $W$ is a semi-direct product $\lie h=\lie g\ltimes_\rho V$ define
the \emph{dual semi-direct product $\lie h^\vee$} by $\lie
h^\vee:=\lie g\ltimes_{\rho^*} V^*$ where $\rho^*$ is the dual
representation
\[
(\rho^*(x)a)(u):=-a(\rho(u)).
\]
Suppose $J$ is totally real with respect to the semi-direct product
$\lie h$.  As noted in (\ref{vanishing components}), the sole
obstruction for the differential of the induced two-form $\omega_J$ on
$\lie h^\vee$ to vanish is due to $(d\omega_J)(x,y,u)$ where $x,y$ are
in $\lie g$ and $u$ is in $V^*$. In the present case,
\begin{align*}
  (d\omega_J)(x,y,u) &= -\left( \omega_J([x,y],u) + \omega_J([y,u],x)
    + \omega_J([u,x],y)\right)\\
    &=-\left( u(J[x,y])-(\rho^*(y)u)(Jx)+(\rho^*(x)u)(Jy)
    \right)\\
  &=  -\left( u(J[x,y]) + u(\rho(y)Jx) - u(\rho(x)Jy)\right)\\
  &=  -u(J[x,y]+ \rho(y)Jx -\rho(x)Jy)\\
  &=  -u(J([x,y] - J\rho(y)Jx +J\rho(x)Jy)).
\end{align*}
By (\ref{eq:5}), we have the following.
\begin{lemma}\label{lem:complex to symplectic}
  Suppose $\lie h$ is a semi-direct product $\lie g\ltimes_\rho V$ and
  let $J$ be a totally real complex structure on $\lie h$.  Then the
  dual semi-direct product $\lie h^\vee$ is Lagrangian with respect to
  the two-form $\omega_J$.  Moreover, $\omega_J$ is symplectic if and
  only if $J$ is integrable.
\end{lemma}
Similarly, if $\omega$ is a non-degenerate two-form on $\lie h$, and
the semi-direct product is Lagrangian, then the Nijenhuis tensor of
the induced complex structure $J_\omega$ is determined by
\begin{align*}
  N_{J_\omega}(x,y) &= [x,y] - [J_\omega(x),J_\omega(y)] +
  J_\omega([x,J_\omega(y)] - [y,J_\omega(x)])\\
  &= [x,y] - [\omega(x),\omega(y)] -
  \omega^{-1}([x,\omega(y)] - [y,\omega(x)])\\
  &= [x,y] - \omega^{-1}(\rho^*(x)\omega(y) - \rho^*(y)\omega(x)).
\end{align*}
By (\ref{closedness}), we have the following.
\begin{lemma}\label{lem:symplectic to complex}
  Suppose $\lie h$ is a semi-direct product $\lie g\ltimes_\rho V$
  with a non-degenerate two-form $\omega$. Suppose that the
  semi-direct product is Lagrangian.  Then the complex structure
  $J_\omega$ on $\lie h^\vee$ is totally real.  Furthermore,
  $J_\omega$ is integrable if and only if $\omega$ is symplectic.
\end{lemma}

In Lemma \ref{lem:complex to symplectic} and Lemma \ref{lem:symplectic
  to complex} (compare~\cite{MR2104679}), we demonstrate the passages
from a totally real complex structure on a semi-direct product to a
symplectic structure on the dual semi-direct product, and from a
symplectic structure on a Lagrangian semi-direct product to a totally
real complex structure on the dual. In the next lemma, we demonstrate
that these two processes reverse each other.

\begin{lemma}\label{inverse}
  Let $\lie h=\lie g\ltimes_\rho V$ be a semi-direct product. Let
  $\phi\colon \lie h\to(\lie h^\vee)^\vee$ be the canonical
  isomorphism defined by the identification $(V^*)^*=V$.  If $\lie h$
  is equipped with a totally real complex structure $J$ then
  $\phi^*(J_{\omega_J})=J$.  Similarly, if $\lie h$ is Lagrangian with
  respect to a symplectic form $\omega$ then
  $\phi^*(\omega_{J_\omega})=\omega$.
\end{lemma}
\bproof The identification $(V^*)^*=V$ is of course the map $u\mapsto
u^{**}$ given by $u^{**}(v^*):=v^*(u)$ for all $u\in V$ and $v^*\in
V^*$.  Setting $\phi(x+u):=x+u^{**}$ this is an isomorphism of Lie
algebras as one may easily check.  Suppose $J$ is totally real.  Then
$J_{\omega_J}$ is also totally real, by Lemmas \ref{lem:complex to
  symplectic} and \ref{lem:symplectic to complex}.  Then the second
statement follows by checking that $J_{\omega_J}x=\phi(Jx)\in
(V^*)^*$.  But
\begin{equation*}
  (J_{\omega_J}x)(u^*)=\omega_J(x)(u^*)=u^*(Jx)=(Jx)^{**}(u^*)=
  \phi(Jx)(u^*).
\end{equation*}
Similarly, if $\lie h$ is Lagrangian with respect to $\omega$ then it
is also Lagrangian with respect to $\omega_{J_\omega}$.  Moreover,
\begin{equation*}
  \omega_{J_\omega}(x,u^{**})=u^{**}(J_\omega x) = (J_\omega x)(u) =
  \omega(x,u).
\end{equation*}
This completes the proof.  
\eproof

\subsection{Special Lagrangian structures}\label{sec:special lag}

A non-degenerate two-form $\omega$ and a complex structure $J$ on a
vector space $W$ are said to be \emph{compatible} if
$\omega(J\xi,J\eta)=\omega(\xi,\eta)$ for all $\xi,\eta$ in $W$.  In
that case $g(\xi,\eta):=\omega(\xi,J\eta)$ is a non-degenerate
symmetric two-tensor on $W$, the induced metric for which $J$ is an
orthogonal transformation: $g(J\xi,J\eta) = g(\xi,\eta)$. If $g$ is
positive-definite we say $(\omega,J)$ is an almost Hermitian pair,
otherwise we say $(\omega, J)$ is almost pseudo-Hermitian.

Suppose that $(\omega, J)$ is an almost pseudo-Hermitian structure on
$W$ and a vector subspace $V$ is totally real with respect $J$.  Then
$V$ is isotropic with respect to $\omega$ if and only if $JV$ is
isotropic.  If $V$ is a totally real subspace, then the splitting
$W=V\oplus JV$ is orthogonal with respect to the induced metric. In
addition, it is clear that $g_{\vert_{JV}}$ is determined by
$g_{\vert_V}$.

Conversely, let $J$ be a complex structure on $W$ and $V$ a totally
real subspace. Any inner product $g$ on $V$ could be extended to $W$
by declaring $g(Jx,Jy)=g(x,y)$ for $x,y\in V$ and $g(x,Jy)=0$. Then
$(\omega,J)$ is an almost pseudo-Hermitian structure on $W$.

If $W$ is a Lie algebra $\lie h=(W,[\cdot,\cdot])$, we say an almost
pseudo-Hermitian pair is pseudo-K\"ahler if $\omega$ is symplectic and
$J$ is integrable.

\begin{definition}\label{special L}
  Let $\lie h$ be a Lie algebra with a semi-direct product structure
  $\lie h=\lie g\ltimes_\rho V$.  Let $(\omega, J)$ be a
  pseudo-K\"ahler structure on $\lie h$. Then $\lie h$ is said to be
  special Lagrangian if $\lie g$ and $V$ are totally real with respect
  to $J$ and Lagrangian with respect to $\omega$. We then also call
  $(\omega, J)$ a special Lagrangian structure on the semi-direct
  product $\lie h$.
\end{definition}

\begin{proposition}\label{prop:special}
  If $(\omega,J)$ is a special Lagrangian structure on a semi-direct
  product $\lie h=\lie g\ltimes_\rho V$, then $(\omega_J, J_\omega)$
  is a special Lagrangian structure on the dual semi-direct product
  $\lie h^\vee=\lie g\ltimes_{\rho^*}V$.
\end{proposition}
\bproof In view of Lemma \ref{lem:complex to symplectic} and Lemma
\ref{lem:symplectic to complex}, the only issue is to verify that for
any $x+u,y+v$ in $\lie g\oplus V^*$,
\[
\omega_J(J_\omega(x+u), J_\omega(y+v))=\omega_J(x+u, y+v).
\]
Given the compatibility of $\omega$ and $J$, the proof is simply a
matter of definitions as given in (\ref{omega-j}) and
(\ref{j-omega}).
\eproof

Other than allowing the metric being pseudo-K\"ahler, Definition
\ref{special L} above is an invariant version of the usual definition
of special Lagrangian structures found in literature on mirror
symmetry if we extend the metric $g$ and the complex structure $J$ to
be left-invariant tensors on the simply connected Lie groups of $\lie
h$ and $\lie g$ (see e.g. \cite{Hitchin}). To illustrate this point,
note that if $\{e^1, \dots, e^n\}$ is an orthonormal basis of $\lie g$
with respect to the (pseudo-)Riemannian metric $g$, set
$u^j=e^j+iJe^j$. Then $\{u^1, \dots, u^n\}$ is a Hermitian basis of
$\lie h$. Then the K\"ahler form $\omega$ is
\[
\omega=i\sum_{j=1}^n u^j\wedge {\overline
  u}^j=i\sum_j(e^j+iJe^j)\wedge(e^j-iJe^j).
\]
The complex volume form is
\begin{eqnarray*}
  \Phi &=& u^1\wedge\dots \wedge u^n =(e^1+iJe^1)\wedge \cdots \wedge
  (e^n+iJe^n)\nonumber\\
  &=& e^1 \wedge\dots \wedge e^ n + i^n Je^1 \wedge\dots \wedge Je^
  n\\
  &&+ \mbox{ terms mixed with both } e^j \mbox{ and } Je^k.
\end{eqnarray*}
When $n$ is odd the real part of $\Phi$ restricts to zero on $V$ and
the imaginary part restricts to a real volume form. Therefore, the
fibers of the quotient map from the Lie group $H$ onto $G$ are special
Lagrangian submanifolds.

\subsection{Flat connections and special Lagrangian
structures}\label{sec:flat again} Suppose $(\omega,J)$ is a special
Lagrangian structure on $\lie h=\lie g\ltimes_\rho V$ and let $g$ be
the induced metric.  Define $\gamma(x):=-J\rho(x)J$. Then $\gamma$ is
a flat torsion-free connection on $\lie g$.  Since $\omega$
is closed
\begin{align*}
  (d\omega)(x,y,Jz) &= - (\omega([x,y],Jz) + \omega([y,Jz],x) +
  \omega([Jz,x],y))\\
  &= - g([x,y],z) + g(\gamma(y)z,x) - g(\gamma(x)z,y)\\
  &= - g([x,y] - \gamma^t(y)x + \gamma^t(x)y , z) = 0,
\end{align*}
and, since $\gamma$ is flat,
\begin{eqnarray*}
  &&- \gamma^t([x,y]) - \gamma^t(x)\gamma^t(y) +
  \gamma^t(y)\gamma^t(x) \\
  &=& - (\gamma([x,y] - \gamma(x)\gamma(y) + \gamma(y)\gamma(x))^t=0.
\end{eqnarray*}
Therefore, $-\gamma^t$ is another flat torsion-free connection.

On the other hand, suppose that $\lie g$ is equipped with a
non-degenerate bilinear form $g$. Let $\gamma$ be a flat torsion-free
connection such that $\gamma':=-\gamma^t$ is also an flat torsion-free
connection. Then, as above the complex structure on $\lie h := \lie
g\ltimes_\gamma\lie g$ given by $J(x,y)=(-y,x)$ is integrable. We
write $x+Jy,~x,y\in\lie g$ for the elements in $\lie h$.  Define
$\omega$ on $\lie h$ by $\omega(x,y)=\omega(Jx,Jy)=0$ and
$\omega(x,Jy) = g(x,y) = - \omega(y,Jx)$ and set $g(Jx,Jy) = g(x,y)$.
Then essentially the same calculation as above shows that $d\omega=0$
by virtue of $\gamma'$ being flat and torsion-free.

\begin{remark}
  Note that we may equally well choose to work with the integrable
  complex structure $J'(x,y)=(-y,x)$ on $\lie h:= \lie
  g\ltimes_{\gamma'}\lie g$ and the associated symplectic form
  $\omega'$.  This is of course precisely the ``mirror image'' of
  $(\lie h,J,\omega)$.  This all amounts to
\end{remark}

\begin{proposition}
  Let $\lie g$ be a Lie algebra with a non-degenerate bilinear form
  $g$. Then there is a two-to-one correspondence between special
  Lagrangian structures on a semi-direct product extending the Lie
  algebra $\lie g$ and flat torsion-free connections $\gamma$ on $\lie
  g$ such that the dual connection $-\gamma^t$ is also flat and
  torsion-free.
\end{proposition}

\section{Canonical isomorphism of DGAs }\label{sec:canonical}

In this section, we consider the relation between $\DGA(\lie h, J)$
and $\DGA(\lie h^\vee, \omega_J)$ when $\lie h$ is a semi-direct
product totally real with respect to a complex structure $J$.

Let $\gamma$ be a flat torsion-free connection on a Lie algebra $\lie
g$. Write $V$ for the associated representation of $\lie g$ on itself
and consider the usual integrable complex structure $J$ on $\lie
h=\lie g\ltimes_\gamma V$. Then $\lie f^1(\lie h,J) = \lie
h^{(1,0)}\oplus\lie h^{*(0,1)}$ where $\lie h^{(1,0)}$ spanned by
$(1-iJ)x$ as $x$ goes through $\lie g$ while $\lie h^{*(0,1)}$ is
generated by $(1-iJ)\alpha$ where $\alpha$ ranges through
$V^*\subset\lie h^*$. Here $J$ acts on $V^*$ by
$(Jv^*)(x+u)=-v^*(Jx+Ju)=-v^*(Jx)$. In particular, $Jv^*\in
\Ann(V)\subset \lie h^*$.

Now set $\lie h^\vee:= \lie g\ltimes_{\gamma^*}V^*$ and define
 $\phi\colon\lie h^\vee_{\mathbb C} \to \lie f^1(\lie h,J)$ as
the tautological map:
\begin{equation}
\phi(x+v^*):=(1-iJ)x+(1-iJ)v^*.
\end{equation}
Recall that the restriction of the Schouten bracket on the space
$\lie f^1(\lie h, J)$ is a Lie bracket.
\begin{lemma} The map
$
  \phi: \lie h_\CC^\vee \to \lie f^1(\lie h, J)
$
   is an isomorphism of Lie algebras.
\end{lemma}
\bproof
  This is a straight-forward check. First, if $u^*$, $v^*$ are in
  $V^*$, then  $\phi(u^*),\phi(v^*)$ are in $\lie h^{*(0,1)}$.
  Therefore,
  $
  [\phi(u^*),\phi(v^*)]=0=\phi([u^*,v^*]).
  $
  If $x,y\in\lie g$, then
  \begin{align*}
    [\phi(x),\phi(y)] &= [(1-iJ)x,(1-iJ)y]= [x,y] -i
    ([x,Jy]+[Jx,y])\\ & =[x,y] -iJ[x,y] = \phi([x,y])
  \end{align*}
  by integrability of $J$.  Finally, take $x\in \lie g$ and $v^*\in
  V^*$.  Then $[\phi(x),\phi(v^*)]\in \lie h^{*(0,1)}$.   With
  $y\in\lie g$ we get
  \begin{eqnarray*}
  &&  [\phi(x),\phi(v^*)]((1+iJ)y)\\
   &=& [(1-iJ)x,(1-iJ)v^*]((1+iJ)y)\\
    &=& - ((1-iJ)v^*)([(1-iJ)x,(1+iJ)y])\\
    &=& - ((1-iJ)v^*)( [x,y] + i([x,Jy]-[Jx,y]) ).
    \end{eqnarray*}
    It is apparent that $v^*([x,y])=0$. In addition, as $[x,Jy]$ is in $V$
    and $Jv^*$ is in $\Ann(V)$, the above is equal to
    \begin{eqnarray*}
    &=& i(Jv^*)([x,y]) - iv^*([x,Jy]-[Jx,y])\\
    &=& -iv^*(J[x,y]+[x,Jy]-[Jx,y])= -2iv^*([x,Jy]).
  \end{eqnarray*}
  While
  \begin{eqnarray*}
   && \phi([x,v^*])((1+iJ)y) = ((1-iJ)[x,v^*])((1+iJ)y)\\
    &=& -i(J[x,v^*])(y) + i[x,v^*](Jy)= 2i[x,v^*](Jy)=-2iv^*([x,Jy]).
  \end{eqnarray*}
\eproof

Recall that the complex structure $J$ on $\lie h$ induces a
symplectic structure $\omega_J$ on $\lie h^\vee$. Then the
contraction map
\[
\omega_J: \lie h^\vee \to (\lie h^\vee)^*
\]
carries the Lie bracket on $\lie h^\vee$ to a Lie bracket $[-\bullet
-]_{\omega_J}$ on $(\lie h^\vee)^*$. Therefore, $\phi\circ
\omega_J^{-1}$ is a Lie algebra isomorphism from $(\lie
h^\vee)_{\mathbb C}^*$ to $\lie f^1$. It induces an isomorphism from
the underlying Gerstenhaber algebra of $\DGA(\lie h^\vee, \omega_J)$
to that of $\DGA(\lie h, J)$. Next we demonstrate that this map is
also an isomorphism of differential graded algebra. i.e.
\begin{equation}\label{differential iso} \phi\circ \omega^{-1}\circ
  d=\dbar \circ \phi\circ \omega^{-1}.
\end{equation}

We do have an isomorphism at hand. Composing $\phi$ with complex
conjugation on $\lie h_{\CC}^\vee$ and $\lie f^1$ respectively, we get
a complex linear map from $\overline{{\lie h}_{\CC}^\vee}$ to
${\overline{\lie f}}^1$.  But the complexified Lie algebra $\lie
h_{\CC}^\vee$ is isomorphic to $\overline{{\lie h}_{\CC}^\vee}$.  This
yields a Lie algebra isomorphism from $\lie h_{\CC}^\vee$ to
${\overline{\lie f}}^1$.  The dual map induces an isomorphism of the
exterior differential algebra generated by the dual vector spaces and
the corresponding pair of Chevalley-Eilenberg differentials. This
isomorphism \it should be \rm the map given in (\ref{differential
  iso}). To see that $\omega$ plays a proper role, we need more
technical details.

As $ \phi^*:(\lie f^1)^*\to \lie (h^\vee_{\mathbb C})^*, $ the
conjugated map is $ \bar\phi^*:(\bar{\lie f^1})^*\to \lie
(\overline{h^\vee_{\mathbb C}})^*.$ In the next calculation we
implicitly identify the isomorphic Lie algebras $(\lie f^1)^*$ with
$\bar{\lie f^1}$ and $\lie h^\vee_{\mathbb C}$ with its conjugate
$\overline{\lie h^\vee_{\mathbb C}}$.  Hence $\bar\phi^*$ is
identified with the map $ \bar\phi^*:{\lie f^1}\to \lie
({h^\vee_{\mathbb C}})^*.$ Then $ \bar\phi^*\phi$ is a map from $\lie
h^\vee_{\mathbb C}$ to $({h^\vee_{\mathbb C}})^*$. According to
\cite[Proposition 11]{CP}, the map $\phi\circ \omega_J^{-1}$ yields an
isomorphism of differential graded algebra as in (\ref{differential
  iso}) if, up to a constant, $ \bar\phi^*\phi$ is equal to the
contract of $\omega_J$. Therefore, we have the following computation.
\begin{eqnarray*}
 && (\bar\phi^*\phi)(x+u^*)(y+v^*)\\
 &=&   (\bar\phi^*)((1-iJ)x+(1-iJ)u^*)(y+v^*)\\
  &=& \overline{\phi^*((1+iJ)x+(1+iJ)u^*)}(y+v^*)\\
  &=& \overline{\phi^*((1+iJ)x+(1+iJ)u^*)(y+ v^*)}\\
  &=& \overline{((1+iJ)x+(1+iJ)u^*)(\phi(y+ v^*))} \\
  &=& \overline{((1+iJ)x+(1+iJ)u^*)((1-iJ)y+ (1-iJ)v^*))}\\
  &=&((1-iJ)x+(1-iJ)u^*)((1+iJ)y+ (1+iJ)v^*))\\
  &=& 2i(u^*(Jy)-v^*(Jx))= 2i\omega_J(x+u^*,y+v^*).
\end{eqnarray*}
This shows that the isomorphism $\phi$ defines a DGA structure on the
de Rham complex of $\lie h^\vee_{\mathbb C}$ isomorphic to the one
defined by $\omega_J$, since the brackets differ only by
multiplication by a constant.  In particular, we have
\begin{theorem}\label{thm:iso}
  $\DGA(\lie h,J)$ and $\DGA(\lie h^\vee,\omega_J)$ are isomorphic.
\end{theorem}
A similar construction and calculation shows that for a Lagrangian symplectic form
$\omega$ on $\lie h=\lie g\ltimes_\rho V$ the associated differential
Gerstenhaber algebras $\DGA(\lie h,\omega)$ and $\DGA(\lie
h^\vee,J_\omega)$ are isomorphic.

\subsection{Examples}\label{sec:examples}

\subsubsection{K\"ahlerian structure on $\mathbb R^{n}\ltimes \mathbb
R^n$.} \label{h1}

Choose $\lie g=\mathbb R^n$ with trivial Lie bracket.  Then a
representation of $\lie g$ on $V=\mathbb R^n$ is given by a linear map
$\rho\colon\mathbb R^n\to\gl(n,\mathbb R)$.  Pick a basis
$e_1,\dots,e_n$ of $\lie g= \mathbb R^n$.  Represent $\lie g$ on
$V=\mathbb R^n$ with basis $v_1,\dots,v_n$ by declaring $\rho(e_i)$ to
be the diagonal matrix with a $\rho_i$ in the $i$-th place of the
diagonal and zero else. Then $\lie h=\lie g\ltimes_\rho V$ is the Lie
algebra given by the structure equations
\begin{gather*}
  dv^i= - \rho_i e^i\bw v^i
\end{gather*}
where $e^i$ and $v^i$ are the dual elements.  By relabeling the
$v_i$ we may suppose $\rho_1,\dots,\rho_p$ to be the non-zero
structure constants for a certain $p\leq n$. Then
\[ \omega_{a} = \sum_{i=1}^n a_i e^i \bw v^i \]
is symplectic for any $n$-tuple $a=(a_1,\dots,a_n)\in \mathbb R^n$
with $a_i\not=0$. Moreover, the complex structure $J$ defined by
\begin{gather*}
  J(e_i)=v_i \quad \mbox{ and } \quad J(v_i)=-e_i
\end{gather*}
is integrable, since
$
  N_J(e_i,e_j) = J\rho(e_i)v_j - J\rho(e_j)v_i = 0.
  $

It is apparent that  $\omega_a$ is of type $(1,1)$ with respect to
$J$. In fact,
\begin{gather*}
  \omega(e_i,Je_i)=\omega(e_i,v_i)=a_i=\omega(v_i,Jv_i),
\end{gather*}
so $\lie h$ is K\"ahler precisely when all $a_i$ are positive (or all
negative).

To identify the ``mirror image'' $\lie h^\vee$ of $\mathbb
R^{n}\ltimes \mathbb R^n$, we first calculate its Lie bracket
$[-,-]^\vee$. Since both the base algebra and the ideal are abelian,
the only non-trivial brackets are contributed by $e_i\in \mathbb
R^n$ and $v^j\in (\mathbb R^n)^*$. As
\begin{eqnarray*}
 && ([e_i,v^j]^\vee)(v_k) =
  (\rho^*(e_i)v^j)(v_k)=-v^j([e_i,v_k])=(dv^j)(e_i,v_k)\\
  &=& \rho_j (e^j\wedge v^j)(e_i,v_k)=\rho_j \delta_i^j\delta_k^j.
\end{eqnarray*}
Therefore,
\[
[e_i,v^j]^\vee =0 \mbox{ when } i\neq j, \mbox{ and }
[e_i,v^i]^\vee= \rho_iv^i \mbox{ for each } i.
\]
In particular, when $\lie h=\mathbb R^{n}\ltimes \mathbb R^n$ is
given by the trivial representation, which has $\rho_i=0$ for all
$i$, then its corresponding dual direct product $\lie h^\vee$ is
again the trivial algebra.

In all cases, the symplectic form and complex structure are
\begin{gather*}
  \omega_J(e_i,v^j)= -v^j(v_i),\quad J_\omega(e_i+v^j) =
  -\omega^{-1}(v^j)+\omega(e_i) = -\frac1{a_j}e_j + a_iv^i.
\end{gather*}

\subsubsection{An example with solvable base.}
\label{solv ex}

Take $\lie g$ to be the solvable Lie algebra with brackets
\begin{equation*}
  [e_1,e_3]=-e_5,\qquad [e_1,e_5]=e_3,\qquad [e_3,e_5]=0.
\end{equation*}
Using $\{e_2, e_4, e_6\}$ as an ordered basis of a vector space $V$,
we represent elements in $\End(V)$ by matrices. We choose an inner
product on $V$ by declaring this basis orthonormal. Then
$\gamma\colon\lie g\to\End(V)$ defined by
\begin{equation*}
  \gamma(e_1)=\left(
    \begin{smallmatrix}
      0&0&0\\0&0&1\\0&-1&0
    \end{smallmatrix}\right),\qquad \gamma(e_3)=0=\gamma(e_5),
\end{equation*}
is a skew-adjoint representation of $\lie g$ on $V$ with respect to
the standard metric. If we consider $V$ as the underlying vector
space of $\lie g$, it also defines a torsion-free left-invariant
connection on $G$. Now, the non-zero brackets on $\lie h=\lie
g\ltimes_\gamma V$ is
\begin{equation*}
  [e_1,e_3]=-e_5,\qquad [e_1,e_4]=-e_6,\qquad [e_1,e_5]=e_3,\qquad
  [e_1,e_6]=e_4.
\end{equation*}
The reader may now verify that the two-form $\omega=e^{1}\wedge
e^2+e^{3}\wedge e^4+e^{5}\wedge e^6$ is symplectic and the complex
structure $Je_{2j-1}=e_{2j}$ is integrable. It is apparent that
$\omega$ is a positive-definite type-(1,1) symplectic form with
respect to $J$. In other words, $(\omega, J)$ is a K\"ahler structure.

To construct $\lie h^\vee$ explicitly, we take the basis $\langle e_1,
e_3, e_5\rangle\oplus \langle e^2, e^4, e^6\rangle$ and identify the
structure equations.
\begin{equation*}
  [e_1,e_3]=-e_5,\qquad [e_1,e^6]=-e^4,\qquad [e_1,e_5]=e_3,\qquad
  [e_1,e^4]=e^6.
\end{equation*}
To construction the corresponding symplectic structure $\omega_J$
and $J_\omega$, we apply (\ref{omega-j}) and (\ref{j-omega}) to find
that
\begin{eqnarray*}
&&\omega_J=e^1\wedge e_2+e^3\wedge e_4+e^5\wedge e_6,\\
&& J_\omega(e_1)=e^2, \quad J_\omega(e_3)=e^4, \quad
J_\omega(e_5)=e^6.
\end{eqnarray*}
We note that both $\lie h$ and $\lie h^\vee$ clearly also may be
represented as semi-direct products of $\mathbb R$ with $\mathbb R^5$
where $\mathbb R^5$ is represented on $\mathbb R^5$ as a line of
transformations skew-symmetric with respect to the standard metric
$g$.  This precisely follows the prescriptions of \cite{MR0425012} to
make $g$ a flat left-invariant metric on $H$.

\section{Nilpotent algebras of dimension at most
six.}\label{sec:nil}
In this section, we tackle two problems when the algebra $\lie h$ is
a nilpotent algebra whose real dimension is at most six.

\begin{problem} Let $\lie h$ be
a  nilpotent algebra. Suppose that it is a semi-direct product $\lie
h=\lie g\ltimes V$ and totally real with respect to a complex
structure $J$. Identify the corresponding algebra $\lie h^\vee=\lie
g\ltimes V^*$ and the associated symplectic structure $\omega_J$.
\end{problem}

In view of Lemma \ref{inverse}, the above problem is equivalent to
finding the associated complex structure on the dual semi-direct
product when one is given a semi-direct product which is also
Lagrangian with respect to a symplectic structure.

The next problem raises a more restrictive issue.

\begin{problem} Let $\lie h$ be
a nilpotent algebra. Suppose that it is a semi-direct product $\lie
h=\lie g\ltimes V$ and it is special Lagrangian with respect to a
pseudo-K\"ahler structure $(J, \omega)$. Identify the corresponding
algebra $\lie h^\vee=\lie g\ltimes V^*$ and the 
associated pseudo-K\"ahler structure $(J^\vee, \omega^\vee)$.
\end{problem}

In view of the example in Section \ref{h1}, when the algebra $\lie
h$ is abelian, the dual semi-direct product $\lie h^\vee$ is again
abelian. The correspondence from $(J, \omega)$ to $(J_\omega,
\omega_J)$ is also given. Therefore, we shall exclude this trivial
case in subsequent computation although we may include it for the
completeness of a statement in a theorem. The first even dimension
in which a non-abelian nilpotent algebra occurs is four.

\subsection{Four-dimensional case}

There are two four-dimensional non-trivial nilpotent algebras
\cite{GK}.  Only one of them is a semi-direct product, namely the
direct sum of a trivial algebra with a three-dimensional Heisenberg
algebra. It happens to be the only one admitting integrable invariant
complex structures \cite[Proposition 2.3]{Salamon}. Up to equivalence,
there exists a basis $\{e_1, e_2, e_3, e_4\}$ on the algebra $\lie h$
such that the structure equation is simply $[e_1, e_2]=-e_3$. The
corresponding complex structure is
\begin{equation}
 \quad J(e_1)=e_2, \quad J(e_3)=e_4, \quad
J(e_2)=-e_1, \quad J(e_4)=-e_3.
\end{equation}
It is integrable. A symplectic form is
\begin{equation}
\omega=e^1\wedge e^4+e^3\wedge e^2.
\end{equation}
Consider the subspaces
\begin{equation}
\lie g:=\langle  e_2, e_4\rangle, \quad V:=\langle e_1, e_3\rangle.
\end{equation}
They determines a semi-direct product $\lie h=\lie g\ltimes_{\Ad} V$.
It is apparent that this semi-direct product is special Lagrangian the
pair $(\omega, J)$ above. One may now work through our theory to
demonstrate that the mirror image of $(\lie h, J, \omega)$ on $(\lie
h^\vee, J_\omega, \omega_J)$ is isomorphic to $(\lie h, J, \omega)$
itself.

\subsection{Algebraic Aspects}\label{sec:algebraic}
In the next few paragraphs, we identify the six-dimensional
nilpotent algebras which is a semi-direct product of a
three-dimensional Lie subalgebra $\lie h$ and an abelian ideal $V$
by identifying equivalent classes of representation of $\lie h$ on
$V$. Once it is done, the construction of the semi-direct product
with the dual representation follows naturally.  We postpone
geometric considerations to the next section.

Since the adjoint action of the nilpotent Lie algebra $\lie g$ on
the abelian ideal $V$ is a nilpotent representation, by Engel
Theorem there exists a basis $\{e_2, e_4, e_6\}$ of $V$ such that
the matrix of any $\ad x$, $x\in \mathfrak{g}$, is strictly lower
triangular.

In our calculation below, we often express the structure equation on
$\lie h=\lie g\ltimes_{\ad}V$ in terms of the C-E differential on
the dual basis $\{e^1, \dots, e^6\}$. In particular we collect
$(de^1, \dots, de^6)$ in an array.  We shall also adopt the
shorthand notation that when $de^1=e^i\wedge e^j+e^\alpha\wedge
e^\beta$, then the first entry in this array is $ij+\alpha\beta$
\cite{Salamon}. To name six-dimensional algebras, we use the
convention developed in \cite{CFGU}.

\subsubsection{Assume that $\lie g$ is abelian.\ \ } There exists a
basis $\{e_1, e_3, e_5\}$ of $\mathfrak{g}$ such that with respect
to the ordered basis $\{e_2, e_4, e_6\}$ for $V$, the adjoint
representation of $\lie g $ on $V$ is given as below
\begin{equation}\label{abel base}
\rho( e_1)= - \left(
\begin{smallmatrix}
0 & 0 & 0 \\
a & 0 & 0 \\
c & b & 0
\end{smallmatrix}
\right), \quad \rho (e_3)= -\left(
\begin{smallmatrix}
0 & 0 & 0 \\
0 & 0 & 0 \\
d & e & 0
\end{smallmatrix}
\right), \quad \rho( e_5)= - \left(
\begin{smallmatrix}
0 & 0 & 0 \\
0 & 0 & 0 \\
f & 0 & 0
\end{smallmatrix}
\right).
\end{equation}
Up to equivalence, we have the following
\begin{equation}\label{h3}
\lie h_3: \quad
 \rho (e_1)= - \left(
\begin{smallmatrix}
0 & 0 & 0 \\
0 & 0 & 0 \\
0 & 1 & 0
\end{smallmatrix}
\right), \quad \rho (e_3)= - \left(
\begin{smallmatrix}
0 & 0 & 0 \\
0 & 0 & 0 \\
1 & 0 & 0
\end{smallmatrix}
\right) , \quad \rho (e_5)=0.
\end{equation}
\begin{equation}\label{h8-a}
\lie h_8: \quad  \rho (e_1)= 0, \quad \rho (e_3)= - \left(
\begin{smallmatrix}
0 & 0 & 0 \\
0 & 0 & 0 \\
1 & 0 & 0
\end{smallmatrix}
\right),  \quad \rho (e_5)=0.
\end{equation}
\begin{equation}\label{h6-a}
\lie h_6: \quad \rho (e_1)= - \left(
\begin{smallmatrix}
0 & 0 & 0 \\
1 & 0 & 0 \\
0 & 0 & 0
\end{smallmatrix}
\right), \quad \rho (e_3)= -\left(
\begin{smallmatrix}
0 & 0 & 0 \\
0 & 0 & 0 \\
1 & 0 & 0
\end{smallmatrix}
\right), \quad \rho (e_5)=0.
\end{equation}
\begin{equation}\label{h17}
\lie h_{17}: \quad \rho (e_1)= - \left(
\begin{smallmatrix}
0 & 0 & 0 \\
1 & 0 & 0 \\
0 & 1 & 0
\end{smallmatrix}
\right), \quad \rho (e_3)=0,  \quad \rho (e_5)=0.
\end{equation}
\begin{equation}\label{h9}
\lie h_9: \quad \rho (e_1)= - \left(
\begin{smallmatrix}
0 & 0 & 0 \\
1 & 0 & 0 \\
0 & 1 & 0
\end{smallmatrix}
\right), \quad \rho (e_3)= -\left(
\begin{smallmatrix}
0 & 0 & 0 \\
0 & 0 & 0 \\
1 & 0 & 0
\end{smallmatrix}
\right), \quad \rho (e_5)=0.
\end{equation}

\subsubsection{Assume that $\lie g$ is non-abelian.\ \ } In this case,
$\lie g$ is a three-dimensional Heisenberg algebra. Thus there exists
a basis $\{e_1, e_3, e_5\}$ of $\lie g$ such that $[ {e_1} \bullet
{e_3} ]=-e_5$, and a basis $\{e_2, e_4, e_6\}$ of $V$ such that the
adjoint representation of $\lie g$ on $V$ is as follows.
\begin{equation}\label{hei base}
\rho( e_1)= -\left(
\begin{smallmatrix}
0 & 0 & 0 \\
a & 0 & 0 \\
c & b & 0
\end{smallmatrix}
\right),  \rho(e_3)= -\left(
\begin{smallmatrix}
0 & 0 & 0 \\
d & 0 & 0 \\
f & e & 0
\end{smallmatrix}
\right),  \rho(e_5)= -\left(
\begin{smallmatrix}
0 & 0 & 0 \\
0 & 0 & 0 \\
bd-ae & 0 & 0
\end{smallmatrix}
\right).
\end{equation}

If $d\neq 0$, by choosing $\{ae_3-de_1, e_3\}$, we have a new set of
$\{e_1, e_3, e_5\}$ such that $d=0$. If $a\neq 0$, we may consider the
new basis $\{e_1, e_3-\frac{d}{a}e_1, e_5\}$ for $\mathfrak{g}$ and
$\{e_2, ae_4+ce_6, e_6\}$ for $V$ and assume $a=1$, $d=0$ and $c=0$.
Then
\begin{equation*}
\lie h=(0,0,0,12, 13, b14-f23+e34+e25).
\end{equation*}
If $e=0$, it is further reduced to
\begin{equation*}
\lie h=(0,0,0,a12, 13, b14-f23).
\end{equation*}
The following becomes easy to verify.
\begin{equation}\label{h6-h}
\lie h_6: \quad \rho (e_1)= -\left(
\begin{smallmatrix}
0 & 0 & 0 \\
1 & 0 & 0 \\
0 & 0 & 0
\end{smallmatrix}
\right), \quad \rho (e_3)=0, \quad \rho (e_5)=0.
\end{equation}
\begin{equation}\label{h7}
\lie h_7: \quad  \rho (e_1)= -\left(
\begin{smallmatrix}
0 & 0 & 0 \\
1 & 0 & 0 \\
0 & 0 & 0
\end{smallmatrix}
\right), \quad \rho (e_3)= -\left(
\begin{smallmatrix}
0 & 0 & 0 \\
0 & 0 & 0 \\
1 & 0 & 0
\end{smallmatrix}
\right), \quad \rho (e_5)=0.
\end{equation}
\begin{equation}\label{h10}
\lie h_{10}: \quad \rho (e_1)= -\left(
\begin{smallmatrix}
0 & 0 & 0 \\
1 & 0 & 0 \\
0 & 1 & 0
\end{smallmatrix}
\right), \quad \rho (e_3)=0, \quad \rho (e_5)=0.
\end{equation}
\begin{equation}\label{h11}
\lie h_{11}: \quad \rho (e_1)= -\left(
\begin{smallmatrix}
0 & 0 & 0 \\
1 & 0 & 0 \\
0 & 1 & 0
\end{smallmatrix}
\right) \quad \rho (e_3)= -\left(
\begin{smallmatrix}
0 & 0 & 0 \\
0 & 0 & 0 \\
1 & 0 & 0
\end{smallmatrix}
\right) \quad \rho (e_5)=0.
\end{equation}
If $e\neq 0$, consider the new basis
\begin{equation*}
\{\ell_1, \dots, \ell_6\}=\{e^2e_1-bee_3-bfe_5, e_2, ee_3+fe_5,
e^2e_4, e^3e_5, e^4e_6\}.
\end{equation*}
Then the structure equations become $(0,0,0,12,13,34+25). $ Taking
the new dual basis $\{e_2-e_3, e_1, e_2+e_3, -e_4+e_5, e_4+e_5,
2e_6\}$, we find this algebra isomorphic to $\mathfrak{h}_{19+}$.
Note that the algebras $\mathfrak{h}_{19+}$ and $\mathfrak{h}_{19-}$
are isomorphic over $\mathbb{C}$ because one has the map
\begin{equation*}
(e_1, e_2, e_3, e_4, e_5, e_6)\mapsto (e_1, e_2, ie_3, e_4, ie_5,
e_6).
\end{equation*}
For future reference, we note that the representation of $\lie g$ on
$V$ is given as below.
\begin{equation}\label{h19}
\lie h_{19}: \quad \rho (e_1)= -\left(
\begin{smallmatrix}
0 & 0 & 0 \\
1 & 0 & 0 \\
0 & b & 0
\end{smallmatrix}
\right) \quad \rho(e_3)= -\left(
\begin{smallmatrix}
0 & 0 & 0 \\
0 & 0 & 0 \\
f & e & 0
\end{smallmatrix}
\right) \quad \rho (e_5)= -\left(
\begin{smallmatrix}
0 & 0 & 0 \\
0 & 0 & 0 \\
-e & 0 & 0
\end{smallmatrix}
\right).
\end{equation}
The last cases are due to $a=d=0$. If it is not already equivalent to
a previous case, they are equivalent to one of the following.
\begin{equation}\label{h4} \lie h_4:
\quad \rho (e_1)= -\left(
\begin{smallmatrix}
0 & 0 & 0 \\
0 & 0 & 0 \\
0 & 1 & 0
\end{smallmatrix}
\right), \quad \rho(e_3)= -\left(
\begin{smallmatrix}
0 & 0 & 0 \\
0 & 0 & 0 \\
1 & 0 & 0
\end{smallmatrix}
\right), \quad \rho (e_5)= 0.
\end{equation}
\begin{equation}\label{h8}
\lie h_8: \quad \rho (e_1)=0, \quad \rho(e_3)=0, \quad \rho (e_5)=
0.
\end{equation}

\subsubsection{The dual semi-direct products. \ \ }
Next, we go on to identify the Lie algebra structure for $\lie
g\ltimes V^*$ for each representation above.  Recall that $\ad=\rho$
has matrix presentations as given in (\ref{abel base}) and (\ref{hei
  base}) depending on whether the algebra $\lie g$ is abelian or not.
The ordered base for $V$ is given by $\{e_2, e_4, e_6\}$. To express
the dual representation $\rho^*$, we shall do it in terms of the
ordered base $\{e^6, e^4, e^2\}$. It amounts to taking the negative of
the ``transpose with respect to the \it opposite \rm diagonal''.
Explicitly, the representation corresponding to (\ref{abel base}) and
(\ref{hei base}) are respectively,
\begin{equation}\label{mirror abel base}
{\rho^*} (e_1)=-
 \left(
\begin{smallmatrix}
0 & 0 & 0 \\
b & 0 & 0 \\
c & a & 0
\end{smallmatrix}
\right), \quad {\rho^*} (e_3)= -\left(
\begin{smallmatrix}
0 & 0 & 0 \\
e & 0 & 0 \\
d & 0 & 0
\end{smallmatrix}
\right), \quad
 {\rho^*} (e_5)=  -\left(
\begin{smallmatrix}
0 & 0 & 0 \\
0 & 0 & 0 \\
f & 0 & 0
\end{smallmatrix}
\right).
\end{equation}

\begin{equation}\label{mirror hei base}
{\rho^*} (e_1)= -\left(
\begin{smallmatrix}
0 & 0 & 0 \\
b & 0 & 0 \\
c & a & 0
\end{smallmatrix}
\right), \quad {\rho^*} (e_3)= -\left(
\begin{smallmatrix}
0 & 0 & 0 \\
e & 0 & 0 \\
f & d & 0
\end{smallmatrix}
\right), \quad {\rho^*} (e_5)= -\left(
\begin{smallmatrix}
0 & 0 & 0 \\
0 & 0 & 0 \\
bd-ae & 0 & 0
\end{smallmatrix}
\right).
\end{equation}
It is now a straight-forward exercise to find the next proposition.
For instance, to find $\lie h^\vee$ when $\lie h$ is $\lie h_3$ as the
semi-direct product of an abelian algebra with an abelian ideal, we
consider the representation (\ref{h3}). By (\ref{mirror abel base}),
the corresponding dual representation is
\[
{\rho^*} (e_1)=
 -\left(
\begin{smallmatrix}
0 & 0 & 0 \\
1 & 0 & 0 \\
0 & 0 & 0
\end{smallmatrix}
\right),  \quad
 {\rho^*} (e_3)=  -\left(
\begin{smallmatrix}
0 & 0 & 0 \\
0 & 0 & 0 \\
1 & 0 & 0
\end{smallmatrix}
\right), \quad {\rho^*} (e_5)= 0.
\]
By (\ref{h6-a}) we find that $\lie h_3^\vee$ is isomorphic to $\lie
h_6$ as a semi-direct product of an abelian algebra with an abelian
ideal. However, if we consider $\lie h_6$ as the semi-direct product
of a Heisenberg algebra and an abelian ideal, we consider the dual
representation of (\ref{h6-h}). Using (\ref{mirror hei base}) and
re-ordering basis on $V^*$, we find that the algebra $\lie h_6^\vee$
is isomorphic to $\lie h_6$. This example reminds us that semi-direct
product structure on a given Lie algebra is not unique, and hence its
dual semi-direct product structure would change accordingly.
Otherwise, elementary consideration such as above provides all
necessary information to complete the next proposition.

\begin{proposition}\label{prop:table of dual} 
  Suppose that $\lie h$ is a six-dimensional nilpotent Lie algebra
  given as a semi-direct product $\lie g\ltimes _{\rho} V$. Then it is
  one of the algebras given in the left-most column of {\rm Table
    (\ref{tab:algebraic mirror}).} Its dual semi-direct product $\lie
  h^\vee:=\lie g \ltimes_{\rho^*}V^*$ is given in the same table as
  checked.
\begin{equation}
\mbox{\small
\begin{tabular}{|l||l|l|l|l|l||l|l|l|l|l|l|l|}
\hline
 $\mathfrak{h}${$\backslash$} $\lie h^\vee$\ &
$\mathfrak{h}_{3}$ & $\mathfrak{h}_{6}$ & $\mathfrak{h}_{8}$ &
$\mathfrak{h} _{9}$ & $\mathfrak{h}_{17}$ & $\mathfrak{h}_{4}$ &
$\mathfrak{h}_{6}$ & $ \mathfrak{h}_{7}$ & $\lie{h}_8$&
$\mathfrak{h}_{10}$ & $\mathfrak{h}_{11}$ & $\mathfrak{h} _{19}$ \\
\hline\hline
$\mathfrak{h}_{3}$ &  & $\checkmark $ &  &  &  &  &  &  &  &  & &  \\
\hline $\mathfrak{h}_{6}$ & $\checkmark $ &  &  &  &  &  &  &  &  &
&  & \\ \hline $\mathfrak{h}_{8}$ &  &  & $\checkmark $ &  &  &  & &
&  &  &  & \\ \hline $\mathfrak{h}_{9}$ &  &  &  & $\checkmark $ & &
&  &  &  &  & &  \\ \hline $\mathfrak{h}_{17}$ &  &  &  &  &
$\checkmark $ &  &  &  &  &  &  & \\ \hline \hline
$\mathfrak{h}_{4}$ &  &  &  &  &  &  &  &  $\checkmark $ &  &  & &  \\
\hline
$\mathfrak{h}_{6}$ &  &  &  &  &  &  & $\checkmark $ &  &  &  &  & \\
\hline $\mathfrak{h}_{7}$ &  &  &  &  &  & $\checkmark $ &  &  &  &
& & \\ \hline $\mathfrak{h}_{8}$ &  &  &  &  &  &  &  &  &
$\checkmark $ &  & &
\\ \hline $\mathfrak{h}_{10}$ &  &  &  &  &  &  &  &  & &
$\checkmark $ &  &
 \\ \hline $\mathfrak{h}_{11}$ &  &  &  &  &  & & &
&  &  & $\checkmark $ &  \\ \hline $\mathfrak{h}_{19}$ &  &  &  &  &
& &  &  &  &  &  & $\checkmark $ \\ 
\hline
\end{tabular}
}
 \label{tab:algebraic mirror}
\end{equation}
\end{proposition}

In {\rm Table (\ref {tab:algebraic mirror})}, the upper left corner is
due to the correspondence between semi-direct products of a
three-dimensional abelian algebra with an abelian ideal. The lower
right corner is due to semi-direct products of a three-dimensional
Heisenberg algebra with an abelian ideal.

\subsection{Geometric Aspects}\label{sec:geometric}

Recall that if $\lie h=\lie g\ltimes V$ admits a totally real
integrable complex structure, then $\lie h^\vee$ is Lagrangian with
respect to a symplectic structure $\omega_J$. As $\lie h_6^\vee=\lie
h_3$ when the base is abelian and $\lie h_3$ does not admit
invariant symplectic structure, $\lie h_6$ as the semi-direct
product of an abelian ideal with an abelian subalgebra would not
admit totally real integrable complex structure. For the same
reason, $\lie h_{19}$ does not admit compatible complex structure.
On the other hand, $\lie h_{17}$ simply would not admit any complex
structure \cite{Salamon}.

When $\lie h_8$ is the semi-direct product of the Heisenberg algebra
with an abelian ideal as given in (\ref{h8}), the integrability of a
compatible integrable complex structure as given in (\ref{eq:5})
implies that the algebra $\lie g$ is abelian. This contradiction
implies that when $\lie h_8$ admits a semi-direct product structure
with a compatible complex structure, then it is the semi-direct
product of an abelian subalgebra and an abelian ideal.

 Therefore, the potential identification from
$(\lie h, J)$ to $(\lie h^\vee, \omega)$ is reduced to the next
table.
\begin{equation}
\mbox{
\begin{tabular}{|l||l|l|l||l|l|l|l|l|l|}
\hline
 $(\mathfrak{h}, J)$  {$\backslash$}  $(\lie h^\vee, \omega)$\ &
 $\mathfrak{h}_{6}$ & $\mathfrak{h}_{8}$ & $\mathfrak{h} _{9}$ &
 $\mathfrak{h}_{4}$ & $\mathfrak{h}_{6}$ & $ \mathfrak{h}_{7}$ &
$\mathfrak{h}_{10}$ & $\mathfrak{h}_{11}$  \\
\hline\hline
$\mathfrak{h}_{3}$   &$\checkmark $ &  &  &  &  &  &  &  \\
\hline   $\mathfrak{h}_{8}$ &    & $\checkmark $ &   & & & & &
\\ \hline $\mathfrak{h}_{9}$ &    &  & $\checkmark $ &  & & & &
\\  \hline \hline
$\mathfrak{h}_{4}$ &  &    &  &  &  &  $\checkmark $ &  &     \\
\hline
$\mathfrak{h}_{6}$ &  &    &  &  & $\checkmark $ &  &  &    \\
\hline $\mathfrak{h}_{7}$ &  &    &  & $\checkmark $ &  &  &  &
\\  \hline $\mathfrak{h}_{10}$ &  &  &    &  &  &  &
$\checkmark $ &
 \\ \hline $\mathfrak{h}_{11}$ &  &    &  & &
&  &  & $\checkmark $   \\  \hline
\end{tabular}
}
\label{tab:geometric mirror}
\end{equation}

\subsubsection{Totally real semi-direct products}\label{real}
Among all the algebras identified in Table (\ref{tab:geometric
mirror}) above, the representation of $\lie g$ on $V$ has the form
\begin{equation}
{\rho} (e_1)= -\left(
\begin{smallmatrix}
0 & 0 & 0 \\
A & 0 & 0 \\
0 & B & 0
\end{smallmatrix}
\right), \quad {\rho} (e_3)= -\left(
\begin{smallmatrix}
0 & 0 & 0 \\
0 & 0 & 0 \\
C & 0 & 0
\end{smallmatrix}
\right), \quad {\rho} (e_5)= 0.
\end{equation}
where $A, B, C$ are respectively zero or one. In addition, $ [e_1,
e_3]=-D e_5$ where $D$ is equal to zero or one, depending on whether
$\lie g$ is abelian or not. In particular, the potentially
non-trivial structure equations are
\[
[e_1, e_2]=-Ae_4, \quad [e_1, e_4]=-Be_6, \quad [e_3, e_2]=-Ce_6,
\quad [e_1, e_3]=-D e_5.
\]

If $J$ is a complex structure such that the semi-direct product is
totally real, then there is a $3\times 3$-matrix $(a_{ij})$ such
that $Je_{2i-1}=\sum_ja_{ij}e_{2j}$. In addition, it is integrable
if and only if $N_J(e_1, e_3)=N_J(e_1, e_5)=N_j(e_3, e_5)=0$. Given
these constraints, one could apply elementary method to identify the
set of complex structures for each set of parameters $(A, B, C, D)$
corresponding to an algebra in the left-most column of Table
(\ref{tab:geometric mirror}). We leave it as exercise.  Instead we
focus on special phenomena.

The first special case is when the algebra is $\lie h_3$ as it does
not admit invariant symplectic structure. However, due to a
classification of complex structure \cite[Proposition 3.4]{Salamon},
up to equivalence, there is a unique complex structure on $\lie
h_3$. With respect to our notations here, it is given by $A=0, B=1,
C=1, D=0$. Making use of the dual representation ordered as in
(\ref{mirror abel base}), and with respect to the bases $\{e_1, e_3,
e_5\}$ on $\lie g$ and $\{e^6, e^4, e^2\}$ on $V^*$, the structure
equation on ${\lie h}_3^\vee=\lie h_6$ is
\[
[e_1, e^6]=-e^4, \quad [e_3, e^6]=-e^2.
\]
Then the 2-form on $\omega_J=e^1\wedge e_2+e^3\wedge e_4+e^5\wedge
e_6$ is a symplectic form on $\lie h_6$.

The second special case is concerned with  $\lie h_6$.
 This
algebra as a semi-direct product is given by  $A=D=1$ and $B=C=0$.
In dual form, the structure equations in the present coordinates are
\[
de^4=e^{12} \quad \mbox{ and }\quad de^5=e^{13}.
\]
It follows that  the constraints for $J$ to be integrable are
\begin{equation}\label{a3}
a_{31}=0, \quad a_{32}=a_{21}, \quad a_{33}=0.
\end{equation}
Therefore,
\begin{equation}\label{Je5}
Je_5=a_{32}e_4, \quad \mbox{ or } \quad a_{32}Je_4=-e_5.
\end{equation}
In particular, $a_{21}=a_{32}\neq 0$.

Let $\omega$ be a symplectic structure on $\lie h_6$ such that $\lie
g$ and $V$ are both Lagrangian.  If $b_{ij}:=\omega(e_i,e_j)$, then
by using that $\omega$ is closed we obtain that $b_{54}=b_{56}=0$
(recall that $e_5\in[\lie h_6,\lie h_6]$ and $e_4,e_6$ belong to the
center) and $b_{52}=b_{43}$.  Let us now assume that $\omega$ and
$J$ are compatible. Then $ \omega(Je_4, Je_3)=\omega(e_4, e_3). $ It
is equivalent to
\[
\omega(a_{32}Je_3, Je_3)=a_{32}\omega(e_4,e_3)=a_{32}b_{43}.
\]
By (\ref{Je5}) above,
\begin{eqnarray*}
&&a_{32}b_{43}=-\omega(e_5, Je_3)\\
&=&-\omega(e_5, a_{21}e_2+a_{22}e_4+a_{23}e_6)=-a_{21}b_{52}-a_{22}b_{54}-a_{23}b_{56}\\
&=&-a_{21}b_{52}=-a_{21}b_{43}=-a_{32}b_{43}.
\end{eqnarray*}
Since $a_{32}\neq 0$, it is possible only when $b_{43}=0$. As
$b_{52}=b_{43}=0$ and $b_{54}=b_{56}=0$, $\omega$ would have been
degenerate. It should that $h_6$ does not admit any special
Lagrangian structure with respect to any semi-direct product
decomposition.

\subsubsection{Family of special Lagrangian algebras}\label{family}
In this paragraph, we  establish the existence of special Lagrangian
structures on the algebras $\lie h_4, \lie h_7, \lie h_9, \lie
h_{10}$ and $\lie h_{11}$. As it turns out, they could be considered
as a family of special Lagrangian structures with ``jumping''
algebraic Lie structures, and hence jumping complex and symplectic
structures.

 We fix a basis $\{
e_1,e_2,e_3,e_4,e_5,e_6\}$ of a real vector space $\lie h$ and
consider also fixed structures $J$ and $\omega$ defined by
\begin{equation}\label{family of j and omega}
Je_{2j-1}=e_{2j}, \qquad Je_{2j}=-e_{2j-1} , \qquad
\omega=e^{16}-e^{25}+e^{34}.
\end{equation}
The 2-form $\omega$ is type (1,1) with respect to $J$, and the
non-degenerate symmetric bilinear form $g(-,-):=\omega(-,J-)$ has
signature $(4,2)$. If $\lie g=\langle e_1,e_3,e_5\rangle$ and
$V=\langle e_2,e_4,e_6\rangle$ then $\lie g$ and $V$ are totally
real with respect to $J$ and maximally isotropic with respect to
$\omega$.

Let $(a,b,c,d)$ be real numbers. For each member of the family of
Lie brackets
\[
(0,0,0,a12,b13,c14+d23),
\]
the corresponding Lie algebra $\lie h$ is the semi-direct product
$\lie h=\lie g\ltimes V$, and the ideal $V$ is abelian.  The
constraint on $\lie h_{a,b}$ so that $\omega$ is closed is
equivalent to $a+b+d=0$. To find the constraints on $\lie h_{a,b}$
so that  $J$ is integral, we choose a basis for the $(1,0)$-forms
with $\omega^j=e^{2j-1}+ie^{2j}$, $1\leq j\leq 3$. As $d\omega^1=0$
and $d\omega^2$ is type $(1,1)$, the sole constraint is due to
$\omega^1\wedge\omega^2\wedge d\omega^3=0$. Given the structure
equations, the  integrability of $J$ is equivalent to $b-c-d=0$.
 It then follows that for any $a,b\in {\mathbb R}$,
\begin{equation}
\Big(\lie h_{a,b}=(0,0,0,a12,b13,(a+2b)14-(a+b)23),J,\omega\Big),
\end{equation}
is a family of special Lagrangian pseudo-K\"ahler  structures on
nilpotent Lie algebras.  Since a non-zero scalar multiple of a Lie
bracket gives rise to just a homothetic change in the metric, we
will restrict ourselves to the curve $\{\lie
h_{a,b}:a^2+b^2=1,\;b\geq 0\}$.  The following isomorphisms can be
checked by using the new basis on the right:
\[
\begin{array}{l}
  \lie h_{\pm 1,0}\simeq \lie h_9, \quad \{ e_2,e_1,e_5,e_3,-e_4,-e_6\}; \\ \\
  \lie h_{0,1}\simeq \lie h_4, \quad \{ e_1,e_3,e_2,\frac12 e_4,e_5,e_6\}; \\ \\
  \lie h_{-\frac{1}{\sqrt{2}},\frac{1}{\sqrt{2}}}\simeq \lie h_{10}, \quad \{ e_1,e_2,e_3,-\tfrac{1}{\sqrt{2}}e_4,
  \tfrac{1}{\sqrt{2}}e_5,-\tfrac{1}{\sqrt{2}}e_6\}; \\ \\
  \lie h_{-\frac{2}{\sqrt{5}},\frac{1}{\sqrt{5}}}\simeq \lie h_7, \quad \{ e_1,e_2,e_3,-\tfrac{2}{\sqrt{5}}e_4,
  \tfrac{1}{\sqrt{5}}e_5,\tfrac{1}{\sqrt{5}}e_6\}; \\ \\
  \lie h_{a,b}\simeq \lie h_{11} \; {\rm if} \; a,b,a+2b,a+b\ne 0,
  \\ \\
   \hspace{2in}  \{ re_1,\tfrac{1}{r}e_2,\tfrac{1}{r}e_3,
  ae_4,be_5,ra(a+2b)e_6\},
\end{array}
\]
where $r=-\left(\tfrac{a+b}{a(a+2b)}\right)^{1/3}$.

Isomorphism classes of such structures translate themselves in this
context as the orbits of the natural action of the group
\[
\Un (2,1)=\{\varphi\in \GL_6(\mathbb R):\varphi
J\varphi^{-1}=J,\;\omega(\varphi\cdot,\varphi\cdot)=\omega\},
\]
on the space of Lie brackets \cite{Lrt}.  Let $\lie h_{a,b}$, $\lie
h_{a',b'}$ be two points in the curve isomorphic as Lie algebras to
$\lie h_{11}$, and assume there exists $\varphi\in\Un (2,1)$ such
that
\[
\varphi.[\cdot,\cdot]_{a,b}:=\varphi[\varphi^{-1}\cdot,\varphi^{-1}\cdot]_{a,b}=[\cdot,\cdot]_{a',b'}.
\]
By using that $\varphi$ must leave invariant the subspaces (which
coincide for both Lie algebras)
\[
\begin{array}{ll}
[\lie h,\lie h]=\langle e_4,e_5,e_6\rangle_{\mathbb R}, & [\lie h,[\lie h,\lie h]]=\langle e_6\rangle_{\mathbb R}, \\ \\
\lie z=\langle e_5,e_6\rangle_{\mathbb R} \; {\mbox (\rm center)}, &
\{ x\in\lie h:[x,\lie h]\subset\lie z\}=\langle
e_3,e_4\rangle_{\mathbb R}
\end{array}
\]
as well as their orthogonal complements relative to both $\omega$
and the metric $g=\omega(\cdot,J\cdot)$, one can easily show that
the matrix of $\varphi$ with respect to the basis $\{
e_1,e_2,...,e_6\}$ is necessarily diagonal, and so it has the form
\[
\varphi=\left[\begin{smallmatrix}
  s &  &  & &  &  \\
   & s &  &  &  &  \\
   &  & \pm 1 &  &  &  \\
   &  &  & \pm 1  &  \\
   &  &  &  & 1/s &  \\
   &  &  &  &  & 1/s \end{smallmatrix}\right],\qquad s\ne 0.
\]
It follows that $a'=\pm\tfrac{1}{s^2}a$, $b'=\pm\tfrac{1}{s^2}b$ and
hence $(a',b')=\pm (a,b)$.

We summarize the results obtained above in the following
proposition.

\begin{proposition}\label{curve} With complex structure $J$ and
symplectic structure $\omega$ given in {\rm (\ref{family of j and
omega})}, with $t\in[0,\pi]$ the function
\[
\Big((0,0,0, (\cos{t})12,
(\sin{t})13,(\cos{t}+2\sin{t})14-(\cos{t}+\sin{t})23),J,\omega\Big)
\]
determines a closed curve in the moduli space of isomorphism classes
of  special Lagrangian pseudo-K\"ahler structures on the space of
six-dimensional nilpotent Lie algebras.  It contains exactly one
structure on each one of {\rm $\lie h_4$ ($t=\pi/2$), $\lie h_{7}$
($t=\arctan{(-\frac12)}$), $\lie h_9$ ($t=0,\pi$)  $\lie h_{10}$
($t=3\pi/4$)}, and the remaining is a continuous family on $\lie
h_{11}$.
\end{proposition}

\begin{theorem}\label{thm:mirror} Suppose that a six-dimensional
  nilpotent algebra admits a special Lagrangian semi-direct product
  structure, then the algebra is isomorphic to one of the following:
  $\lie h_1, \lie h_4, \lie h_7, \lie h_8, \lie h_9, \lie h_{10}, \lie
  h_{11}$. Their weak mirror images are respectively given below.
\begin{center}
\begin{tabular}{|l||l|l|l|l|l|l|l|l|}
  \hline $(\mathfrak{h}, J, \omega)$ & $\lie h_1$& $\lie h_4$ & $\lie
  h_7$ & $\lie h_8$ & $\lie h_9$  & $\lie h_{10}$ & $\lie h_{11}$
  \\
  \hline
  $({\lie h}^\vee, {\omega}_J, {J}_\omega)$ &
  $\mathfrak{h}_{1}$&
  $\mathfrak{h}_{7}$ &  $ \mathfrak{h}_{4}$ & $\mathfrak{h}_{8}$ & $\mathfrak{h} _{9}$  &
  $\mathfrak{h}_{10}$ & $\mathfrak{h}_{11}$  \\
  \hline
\end{tabular}
\end{center}
In addition,
\begin{itemize}
\item All special Lagrangian semi-direct product structure
  $(J,\omega)$ on $\lie h_1, \lie h_8$, $\lie h_9$ and $ \lie h_{10}$
  are self-mirror. It means that there is an isomorphism $\lie
  h_{\ell}\cong \lie h_\ell^\vee$ and there are quasi-isomorphisms.
  \[\DGA({\lie h}_\ell^\vee, \omega_J)\approx\DGA(\lie h_\ell,
  \omega)\approx\DGA({\lie h}_\ell^\vee, J_\omega)\approx \DGA(\lie
  h_\ell, J).\]
\item The algebras $\lie h_7$ and $\lie h_4$ form a mirror pair of
  special Lagrangian semi-direct products.
\item The mirror of a special Lagrangian semi-direct structure on
  $\lie h_{11}$ is again a special Lagrangian semi-direct structure on
  $\lie h_{11}$, but the two pseudo-K\"ahler structures are not
  equivalent.
\end{itemize}
\end{theorem}
\bproof The statements on $\lie h_1$ are trivial.

As noted in Section \ref{real}, the only algebras admitting special
Lagrangian structures are given in Table (\ref{tab:geometric
mirror}).  $\lie h_3$ is eliminated because it does not admit
invariant symplectic structure.

The construction in this section and the paragraph concerned with a
construction on $\lie h_6$ at the end of  Section \ref{real}
demonstrate the ``existence'' of special Lagrangian structures on the
named algebras.

The isomorphisms between $\DGA(\lie h, J)$ and $\DGA(\lie h^\vee,
\omega_J)$ are given by Theorem \ref{thm:iso} and the identification
of algebraic mirrors ${\lie h}^\vee$ as given in Table
(\ref{tab:geometric mirror}).  The validity of the claim on $\lie
h_{11}$ is due to an analysis in \cite[Section 5.5]{CP}.\eproof

\noindent{\bf Acknowledgment\ }\  Y.~S.~Poon would like to thank
Conan N.~C.~Leung for several stimulating conversations.  R.~Cleyton
was supported by the Junior Research Group ``Special Geometries in
Mathematical Physics'' of the Volkswagen Foundation and the SFB 647
``Space--Time--Matter'' of the DFG.  R. Cleyton thanks the University
of California Riverside for hospitality and support during his visit
there in 2007.  Thanks are also due to Jens Heber and Isabel Dotti for
pointing us to some references we might otherwise have missed.

\end{document}